\documentclass[aoas,preprint]{imsart}
\RequirePackage[OT1]{fontenc}
\RequirePackage{amsthm,amsmath}
\RequirePackage[numbers]{natbib}
\RequirePackage[colorlinks,citecolor=blue,urlcolor=blue]{hyperref}
\usepackage{amssymb,makeidx,latexsym,rotate,epsf,amssymb,graphicx,color, subfigure, url}

\newtheorem{defi}{Definition}
\newtheorem{rem}{Remark}

\arxiv{arXiv:0000.0000}

\startlocaldefs
\numberwithin{equation}{section}
\theoremstyle{plain}

\endlocaldefs

\begin{document}

\begin{frontmatter}
\title{Letter to the Editor}
\runtitle{Letter to the Editor}
%\thankstext{T1}{Footnote to the title with the ``thankstext'' command.}

\begin{aug}
\author{\fnms{Milan} \snm{Stehl\'\i k}\thanksref{m1,m2}\ead[label=e1]{E-Mail: mlnstehlik@gmail.com}}
\and
\author{\fnms{Philipp} \snm{Hermann}\thanksref{m2}\ead[label = e2]{E-Mail: philipp.hermann@jku.at}},

% \ead[label=u1,url]{http://www.foo.com}}

%\thankstext{t1}{Some comment}
%\thankstext{t2}{First supporter of the project}
%\thankstext{t3}{Second supporter of the project}
\runauthor{M. Stehl\'\i k and P. Hermann} % L. St\v relev

\affiliation{Universidad T\'ecnica Federico Santa Mar\'{\i}a\thanksmark{m1} and  Johannes Kepler University Linz\thanksmark{m2}}

\address{Departamento de  Matem\'atica \\ Universidad T\'ecnica Federico Santa Mar\'{\i}a\\ Casilla V 110, Valpara\'iso,  Chile \\
% \printead{e2}\\
\printead*{e1}}

\address{Department of Applied Statistics \\Johannes Kepler University Linz\\ Altenbergerstrasse 69, 4040 Linz, Austria\\
% \printead{e1}\\
\printead*{e2}}
\end{aug}

%\begin{abstract}
%ADD ABSTRACT
%\end{abstract}

%\begin{keyword} \kwd{ADD KEYWORDS}  \end{keyword}

\end{frontmatter}

\section{Dedication} \emph{This Letter is dedicated to the 50th anniversary of unexpected death of  Samuel Stanley Wilks. To exact distribution of His,  Wilks's, statistics first author devoted his "Lambert W research" in 2000-2003. }

\section{Introduction}

In 1938, Samuel Stanley Wilks  proved the $\chi^2$-asymptotics of $-2\ln\Lambda$, where $\Lambda$ is the likelihood ratio statistics in regular exponential family (see Wilks S.S. (1938) \cite{Wilks}).
But how does the exact CDF of $-2\ln\Lambda$ look like?  {Stehl\'\i{k} M. (2003)} \cite{Stehlik03}, derived the exact cumulative distribution function of $-2\ln\Lambda$ and decomposition of \emph{Kullback-Leibler}-divergence (\emph{I}-divergence) in the sense of P\'{a}zman A. (1993) \cite{Pazman1}, by substantial usage of Lambert W function, firstly introduced by Johann Heinrich Lambert in 1758 (see \cite{Lambert}), a contemporary of Euler.
The paper by Goerg G. M. (2011) \cite{Georg11}, \emph{"Lambert W random variables-a new family of generalized skewed distributions with applications to risk estimation"}, introduced a class of so called Lambert W$\times$ F random variables,
\begin{eqnarray}
Y_\gamma:=X\exp(\gamma X),\label{LWF}
\end{eqnarray}
where $\gamma\in R$ is skewness parameter and $X$ is continuous random variable.
{Stehl\'\i{k} M. (2003)} \cite{Stehlik03}  derived the exact distribution of Wilks  statistics $-2\ln\Lambda$ to test for  the scale hypothesis
$H_0:\theta=\theta_0 \ \mbox{versus}\ H_1:\theta\ne\theta_0$
in the regular Gamma family and  proven that Wilks  statistics $-2\ln\Lambda$ is a function of a random variable
 \begin{eqnarray}
 G_u(X)=X-u\ln(X),
 \end{eqnarray}  where  $X$ is random variable from exponential family.
Here, notice that
\begin{eqnarray}
-u\ln Y_{-\frac1u}=G_{u}(X), \ \mbox{for } X>0 \label{Rovnako}
\end{eqnarray}
where $Y_\gamma$ of Goerg G. M. (2011) \cite{Georg11} is defined by (\ref{LWF}) and $\gamma={-\frac1u}.$
The statistical application of the class (\ref{LWF}) and "Lambert W function" is intrinsically related to \emph{I}-divergence decompositions and the importance they play in statistical inference. {Stehl\'\i{k} M. (2003)} \cite{Stehlik03} derived that  Kullback-Leibler divergence in the sense of P\'{a}zman (1993) \cite{Pazman1} has the form \begin{eqnarray}\label{Div}
I_{N}(y,\theta)= \sum_{i=1}^N \{G_{u}(\theta y_i)-G_{u}(u)\},
\end{eqnarray}
$y=(y_1,...,y_N).$  Notice, that $I_{1}(X,1)=G_u(X)-G_u(u),$ is the "basic" information of LR test, based on  just a single random variable $X$,  directly relating nonlinearly transformed $Y_\gamma$ of Goerg G. M. (2011) \cite{Georg11} to $G_u(X)$ of {Stehl\'\i{k} M. (2003)} \cite{Stehlik03} (see (\ref{Rovnako})).
In {Stehl\'\i{k} M. (2006)} \cite{Stehlik06} and  {Stehl\'\i{k} M. (2008)} \cite{Stehlik08} extension of results  to Weibull and generalized Gamma distributions (Ggds) was made. Considered  Ggd  covers for various choices of parameters  of one-sided normal, $\chi^2_n,$   Weibull and in the limit  a log-normal distribution.
The LW function approach based on $G_{u}(X)$ transformation was used for exact inference for Pareto heavy tailed distribution in Stehl\'\i k M. et al. (2010) \cite{Stehlik08CS}.
The LW function approach and $G_{u}(X)$  transformation was used fundamentally in Balakrishnan  and Stehl\'\i k (2008) \cite{Bala} for  extension of results  also to cases of Type I and Type II censored samples and missing data.

In this letter we discuss several important   methodological and practical aspects of Lambert W variable.
According Goerg G. M. (2011) \cite{Georg11} the Lambert W framework is a new generalized way to analyze skewed, heavy-tailed data.
In the next two sections we discuss both, heavy-tails and skewness perspectives of this Lambert W framework.

In the next section "Heavy Tails: On three regimes of \texttt{IGMM}-algorithm", based on heavy-tailedness we define three Regimes of  Goerg G. M. (2011)'s Algorithm 3, and its implementation \texttt{IGMM} in R-package \texttt{LambertW}. However, current implementation of  algorithm 3 cannot work in all three Regimes. In Regime III, where no moments of financial data exist, we show that \texttt{IGMM} is not working. Based on simple graphical method  we give a practical guidelines how to discriminate between regimes. Also we introduce a robust tests for normality against heavy tails to enable a formal statistical procedure for the better linking  of a given data to Regimes.
The introduced methodology is illustrated on LATAM data, used by Goerg G. M. (2011). Also suggestion for correction of Algorithm 3 in Regime III is provided.

In the section "Skewness: On asset Returns and t-distribution" we discussed difficulties with symmetrization of data, based on transformation introduced by Goerg G. M. (2011) \cite{Georg11}.

\section{Heavy Tails: On three regimes of \texttt{IGMM}-algorithm}
In this section we describe three regimes of iterative method of moments introduced by Goerg G. M. (2011) \cite{Georg11} (\texttt{IGMM}-Method). The description is  based on approximations by random walk, respectively to heavy-tailedness of input variable $U.$
Such a description is important, in particular for applicability of Algorithm 3 to any financial data, e.g. LATAM data.
The three regimes are defined as follows:
\begin{enumerate}

\item \textbf{Regime I:} distributions $U$ with finite mean and finite variance (here belongs e.g.  student-$t_\nu$-distribution with $\nu>2$)
\item \textbf{Regime II:} distributions $U$ with finite mean but infinite variance (here belongs e.g.  student-$t_\nu$-distribution with $1<\nu \leq 2$)
\item \textbf{Regime III:} distributions $U$ with $E(|U|) = \infty$ and infinite variance (here belongs e.g.  student-$t_\nu$-distribution with $0<\nu\leq 1$). \end{enumerate}

We are showing that algorithm which works in Regime I (because of Strong-Law of Large Numbers) cannot work well in Regime III, since statistical learning in Regime I is related to arithmetic mean, whereas in Regime III to harmonic mean (see Beran, Schell, and Stehl\'\i k (2014) \cite{Beran}).
 Before any further  methodological discussion we provide illustration
of computation with  \texttt{IGMM}-method for the three regimes.
Since in subexponential family Pareto tail is well fitting to student-$t_\nu$ (used also in Goerg G. M. (2011) \cite{Georg11}), we simulate samples of 1000 observations of student-$t_\nu$-distribution with $\nu=1,1.5,5$ degrees of freedom, to represent all three regimes. The same sample size has been used for the computations on the basis of Pareto distribution.
In these regimes we study sensitivity of parameter estimation of $\mu$, $\sigma$ and $\gamma$ of implemented function \texttt{IGMM}.  The procedure for this sensitivity check is conducted as follows:
\begin{enumerate} \item Simulating a sample $U$ from student or Pareto distribution for all three regimes
\item Transformation of  $Y=U\cdot\exp(\gamma\cdot U)\sigma+\mu$ for all samples U
\item Estimation of parameters $(\mu, \sigma, \gamma)$ for transformed samples by usage of \texttt{IGMM}.
\item Repeat steps 2.-3.  for a different values of $\gamma$. \end{enumerate}

The calculated differences between the true values and their estimators are shown in Table \ref{tab:est}. Higher degrees of freedom lead to better approximations of the parameters.  Deviations are higher for increasing  $\gamma$ and lower degrees of freedom. Due to increasing deviations for smaller $\nu$ it can be assumed that \texttt{IGMM}-method works acceptably  for student-\emph{t}-distribution of Regime I, deviations are larger for Regime II and astronomical deviations are received in Regime III.
The similar results are obtained for Pareto distribution. Astronomical deviations of estimation with heavy-tailed distributions $\alpha = 1, 1.5$ are similar to those of student distributions of Regime II and III.

\begin{table}[h]
\centering
\caption{Estimation of parameters with \texttt{IGMM()} for $U$ having student-$t_\nu,$ or Pareto-$\alpha$ distribution}
\label{tab:est}
\begin{tabular}{lccc|lccc}
  \hline
	\multicolumn{4}{c|}{\textbf{Student  $t_\nu$-distribution}} & \multicolumn{4}{c}{\textbf{Pareto-$\alpha$ distribution}} \\\hline
 $\nu$ & $\mu-\hat\mu$ &  $\gamma-\hat\gamma$ & $\frac{\sigma}{\hat\sigma}$ & $\alpha$ & $\mu-\hat\mu$ &  $\gamma-\hat\gamma$ & $\frac{\sigma}{\hat\sigma}$\\
  \hline
  & $\mu = 0.2$ &  $\gamma =0.1$ & $\sigma =1.5$ &&  $\mu = 0.2$ & $\gamma =0.1$ & $\sigma =1.5$   \\ \hline
   5 & 0.0201 & 0.0182 & 1.2535 & 5 & 1.9479 &   0.3433& 0.2331 \\
   1.5 & 0.6061 & 0.3993 & 5.5353 & 1.5 & -5.25$\cdot10^{6}$ &  11.53 & 4.29$\cdot10^{8}$ \\
   1 & {-1.51$\cdot10^{10}$} & 11.6649 & 0.0000  & 1 & 5.24$\cdot10^{24}$ & 0.1504&  1.12$\cdot10^{26}$  \\ \hline
  & $\mu = 0.2$  & $\gamma =0.3$ & $\sigma =1.5$& &$\mu = 0.2$ & $\gamma =0.2$& $\sigma =1.5$  \\ \hline
   5 & 0.0449 & 0.0531 & 1.2054 & 5 & 2.1923& -0.2561& 2.9713\\
   1.5 & 1.58$\cdot10^{12}$ &  -0.0494 & 3.36$\cdot10^{13}$ & 1.5 & 4.49$\cdot10^{11}$ & 0.0501 & 9.54$\cdot10^{12}$\\
   1 & {2.84$\cdot10^{34}$} &  -0.0496 & 6.04$\cdot10^{35}$ & 1 & 7.74$\cdot10^{52}$& 0.0504 & 1.64$\cdot10^{54}$\\ \hline
  & $\mu = 0.2$ & $\gamma =0.5$& $\sigma =1.5$  & & $\mu = 0.2$ & $\gamma =0.25$& $\sigma =1.5$ \\ \hline
   5& 0.1151 & 0.0445 & 1.2027 & 5 & 2.3260& 0.2121 & 0.3343\\
   1.5 & 1.29$\cdot10^{23}$ &  -0.2494 & 2.74$\cdot10^{24}$ & 1.5 & 3.71$\cdot10^{15}$ & 0.0004 & 7.89$\cdot10^{16}$\\
   1 & {8.99$\cdot10^{59}$} &  -0.2497 & 1.91$\cdot10^{61}$& 1& 9.38$\cdot10^{66}$& 0.0005& 1.99$\cdot10^{68}$ \\\hline
\end{tabular}
\end{table}

The astronomical discrepancies in Regime III (the case where no finite expectation exists), i.e. $0<\nu\leq 1$ are theoretically  explained by law of large numbers. Goerg G. M. (2011) \cite{Georg11} has used in his algorithm \texttt{IGMM} intuitively scaled score function, i.e.  $\sigma S_\mu(X),$ of the normal distribution
\begin{equation}
S_\mu(X):=(X-\mu)/\sigma^2, \label{normalscore}
\end{equation}
where mean $\mu$ is taken as a parameter of interest and $\sigma$ is nuisance. Such an algorithm is working  when both mean and  variance are finite, i.e. for $\nu>2.$ However, when only  mean is finite ($1<\nu \leq 2$), the effect of nuisance $\sigma$ is well visible (see Table \ref{tab:est}).
In the case of heavy tailed student ($0<\nu\leq 1$), where both  mean and variance are infinite, the error converges in probability to infinity.
This can be obtained by usage of e.g. Kolmogorov's Strong Law of Large Numbers (LLN) (see e.g. \cite{Sung}). If we have a sample from distribution with infinite mean (e.g.  $t_1$), i.e. Lebesgue integral $\int_R |x|dF(x)=\infty,$ then $\frac1n \sum_{i=1}^n X_i$ will have a finite limit for $n\to \infty$ with probability zero. Such random walk is introduced in step 8 (among others) of Algorithm 3 of Goerg G. M. (2011) \cite{Georg11}, where sample mean and sample deviation updates scale and location parameters.
Therefore we shall expect astronomical numbers in both differences of location parameters $\mu-\hat\mu$ and ratios of scales $\sigma/\hat\sigma$ with probability 1 (see e.g. rows $\nu=1$ for Student $t_\nu$ and rows $\alpha=1$ for Pareto($\alpha$) in Table \ref{tab:est}). Such a divergence is not avoided by step 4, namely $||\tau^{(k)}-\tau^{(k-1)}||>tol.$, in Algorithm 3 of Goerg G. M. (2011) \cite{Georg11}.

Random walk of normal scores (\ref{normalscore}) is the reason for this behavior and it explains  the astronomical errors of magnitude
$10^{59}$ for $t_1$ distribution. Indeed, especially in financial returns (like Asset returns, discussed in section 7.2 of Goerg G. M. (2011) \cite{Georg11}) we shall expect heavy tailed data. Naturally following questions arise: What should be done in such cases? Can we define some procedure how to check when we can apply \texttt{IGMM}?
The answer to this questions is given in the next Sections 3.1 and 3.2.

\subsection{Robust testing for normality against Pareto tail}

First, we shall  testify for the range of Pareto tail parameter $\alpha$ against light tailed normal distribution: for this purpose we need to apply a test for normality against Pareto tails. A consistent and robust test developed as a robust version of Jarque-Bera (JB) test based on the location functional is given by Stehl\'\i k et al (2012) \cite{StehlikCS}. This procedure recognizes in which regime we have our data.
The developed test also works  for arbitrary sample size, which is very practical for financial applications.
For a specific alternatives, also robustified directed Lin-Mudholkar tests (see {Stehl\'\i{k}, Thulin and  St\v relec (2014) \cite{StehlikCIS}) have a good trade-off between power and robustness.
Before using  such a test one shall check for homogeneity in Pareto tail within our financial time series (this is a practical problem, because tail parameters typically varies  during the series). Such testing  procedure is developed, jointly  with likelihood ratio test for simple hypothesis of the Pareto tail $\alpha=\alpha_0$ in Stehl\'\i k et al (2010) \cite{Stehlik08CS} by a  substantial usage of Lambert W-random variable. We applied robustified JB test of Stehl\'\i k et al (2012) \cite{StehlikCS} for simulated data from Regime II and III, and received p-values 0. Thus it is not recommended to apply \texttt{IGMM} to these two regimes.

To explain this fact from the point of view of finance,
we shall realize that LATAM data contains daily log-returns (in percent) of an equity fund investing in Latin America (LATAM) from January 1, 2002 until May 31, 2007. Emerging markets in Latin America
(see e.g. \cite{RamazanFaruk}) can have different properties on left and right tails. It was shown that e.g. Argentina and Brazil have higher estimates of the right tail index than of the left tail index. Therefore, high positive returns are more likely than similar losses in these growing economies. In 2004 it was observed that positive stock return distribution in e.g. Brazil may not have a finite second moment since the estimated extreme value index was around 0.5.
There is (even from 1988) an empirical evidence of non-existence of first moment (see page 9 of \cite{Akgiray}).
Another increase of heavy tailedness of the right tail has been introduced in the years 2004-2015, where high-frequency trading became more present in Latin America.
Analogously, in electric markets less credibility has been given to analysis using empirical means, often quickly replaced by median based techniques (see e.g. \cite{Kim11}).

\subsection{A graphical screening between regimes of \texttt{IGMM}}

In the following, we show how the three Regimes of \texttt{IGMM} can be recognized based on t-Hill plots.
t-Hill estimator is a robust but consistent Pareto tail estimator introduced in Fabi\'an and {Stehl\'\i{k} (2009) \cite{FSt09}
and its consistency for iid sample was proven in Stehl\'\i k et al (2012) \cite{StehlikCIS}, whereas for dependent data in Jordanova, Du\v sek and Stehl\'\i k (2013) \cite{Chemo1}. We base our regimes discrimination on robust t-Hill, so that regime boundaries are not influenced by possible outliers. However, to decrease variability (and  increase efficiency) of specification of type of regime for a given data, we use flexible Harmonic mean estimator introduced in Beran, Schell and Stehl\'\i k (2014) \cite{Beran}.

To recall the Harmonic mean estimator, the next definition follows.
\begin{defi}
We suppose that $\mathbf{X}_1, \mathbf{X}_2, ..., \mathbf{X}_n$ are possibly dependent copies of $\mathbf{X}$ with d.f. $F$, upper order statistics
$$\mathbf{X}_{(1,n)} \leq \mathbf{X}_{(2,n)} \leq ... \leq \mathbf{X}_{(n,n)}.$$ Let us  denote ${\rm RV}_a$   the class of regularly varying functions at infinity, with an index of regular variation equal to $a\in {\Re}$, i.e.\  positive measurable functions $g(\cdot)$ such that for all $x>0$, $g(tx)/g(t)\rightarrow x^a$, as $t\rightarrow\infty$.
\begin{equation}\label{RVcond}
    1 - F \in {\rm RV}_{-\alpha}, \quad \alpha > 0.
\end{equation}
Harmonic Moment tail Index Estimator has the form  $$H_{k, n}^*(\beta) =\frac{1}{\widehat{\alpha}_{k, n}(\beta)}=
\frac{1}{\beta-1}\left \{\left[{\frac{1}{k} \sum\limits^k_{j=1}
\left ( \frac{X_{n-k,n}}{X_{n-j+1,n}}\right )^{\beta-1}}\right]^{-1}-1\right \},$$
where $\beta > 0$ is tuning parameter.

\end{defi}

For $\beta= 2$ we obtain t-Hill, for $\beta= 1$ we have Hill estimator (see Hill (1975) \cite{Hill1975}).
The tuning parameter $\beta$ is regulating the trade-off
between efficiency and robustness. For $\beta > 1$ the effect of large contaminations
is bounded, since the Harmonic Moment Tail Index Estimator benefits
from the properties of the harmonic mean. However, a larger value of $\beta$ also
implies an increased variance. For $\beta < 1$ the Harmonic Moment Tail Index
Estimator also has a higher variance than Hill's estimator.

\begin{rem}
\textbf{Remark on VAR for  LATAM returns}

As the second example, Goerg G. M. (2011) \cite{Georg11}  reexamines the \texttt{LATAM} returns.  He assures that "a comparison of risk estimators (Value at Risk, VAR) demonstrates the suitability of the Lambert W $\times$ F distributions to model
financial data." From the perspective of minimal mean square error, a Mean-of-order-\emph{p} (MOP) class of VAR estimators
can have a mean square error smaller than that of classical extreme value index (EVI) estimators, not only around optimal levels,  but for other levels too (see Gomes, Brilhante and Pestana (2014) \cite{VARMOP}). MOP EVI-estimator $H_{k,n}^{(p)}$ was introduced in Brilhante, Gomes  and Pestana (2012) \cite{BGP2012}.
Note that if we consider a generalization (motivated by robustness) to $p<0$ of the MOP functionals $H_{k,n}^{(p)}$,  we get
the t-Hill estimator $ H_{k,n}^*(2)=H_{k,n}^{(-1)}$. This is a VAR-justification of why to use t-Hill estimator for specification of boundaries of the Regimes.
Such setup is also of interest for BASEL II (and higher) initiative in banking and audit.
\end{rem}

Let $n$ be fixed as sample size.  Analogously to the Hill plot we consider
the set of points with coordinates
$$\left( k,\, \frac{1}{\widehat{\alpha}_{k, n}(\beta)}\right), \quad k \in \{1, 2, ..., n\}.$$
Further on we call this plot "modified Hill plot". Our graphical procedure is illustrated on discrimination between $t_1,t_{1.5}$ and $t_5$ in Figure \ref{fig:IGMMconv}.

\begin{figure}[h]
%\vspace{6pc}
		\subfigure[Comparison of three t-Hill plots lines]{\includegraphics[width = 0.48\textwidth]{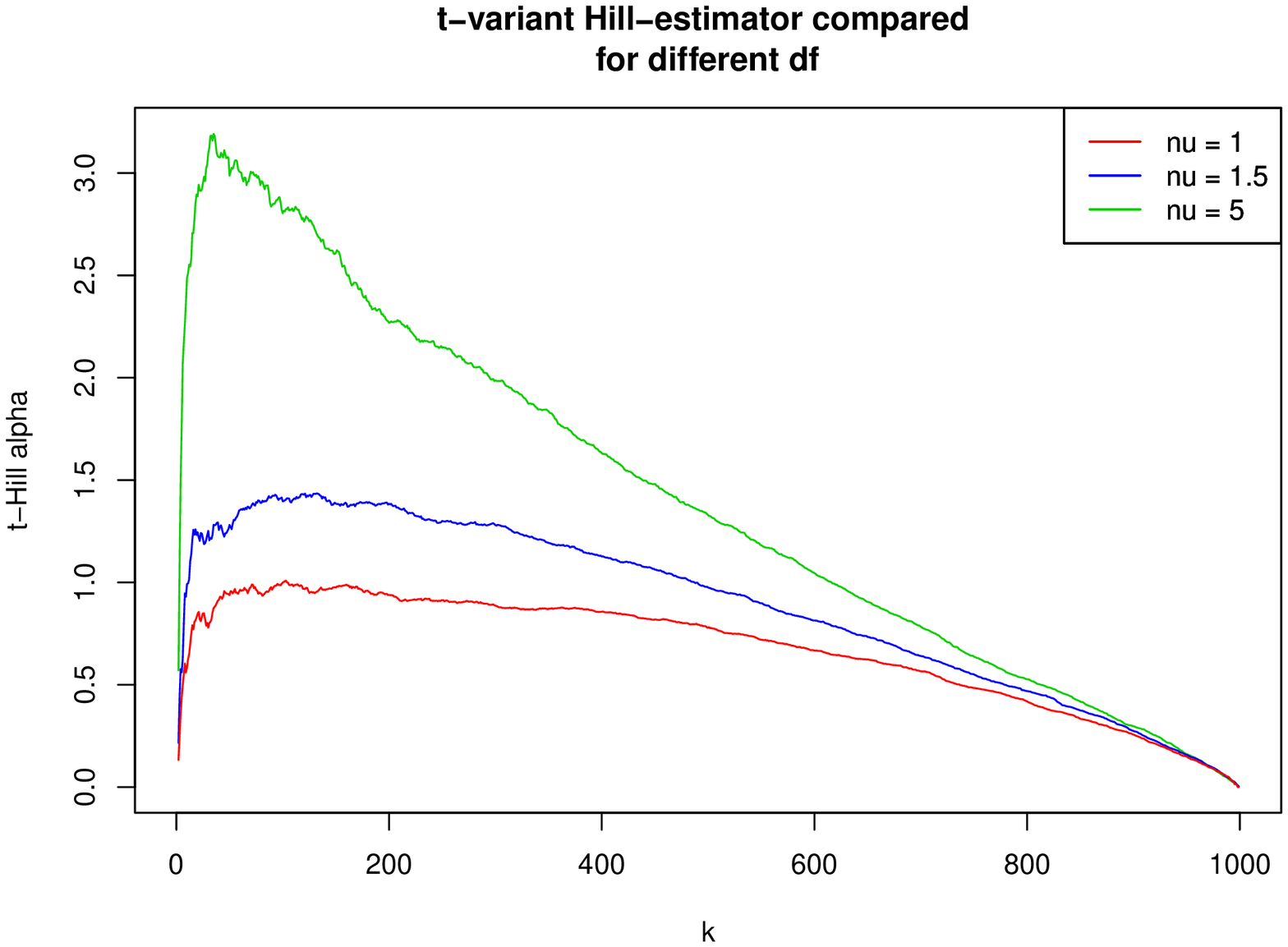}}
		\subfigure[Convergence region of  t-Hill plot: best estimation and distinguishing of 3 regimes]{\includegraphics[width = 0.48\textwidth]{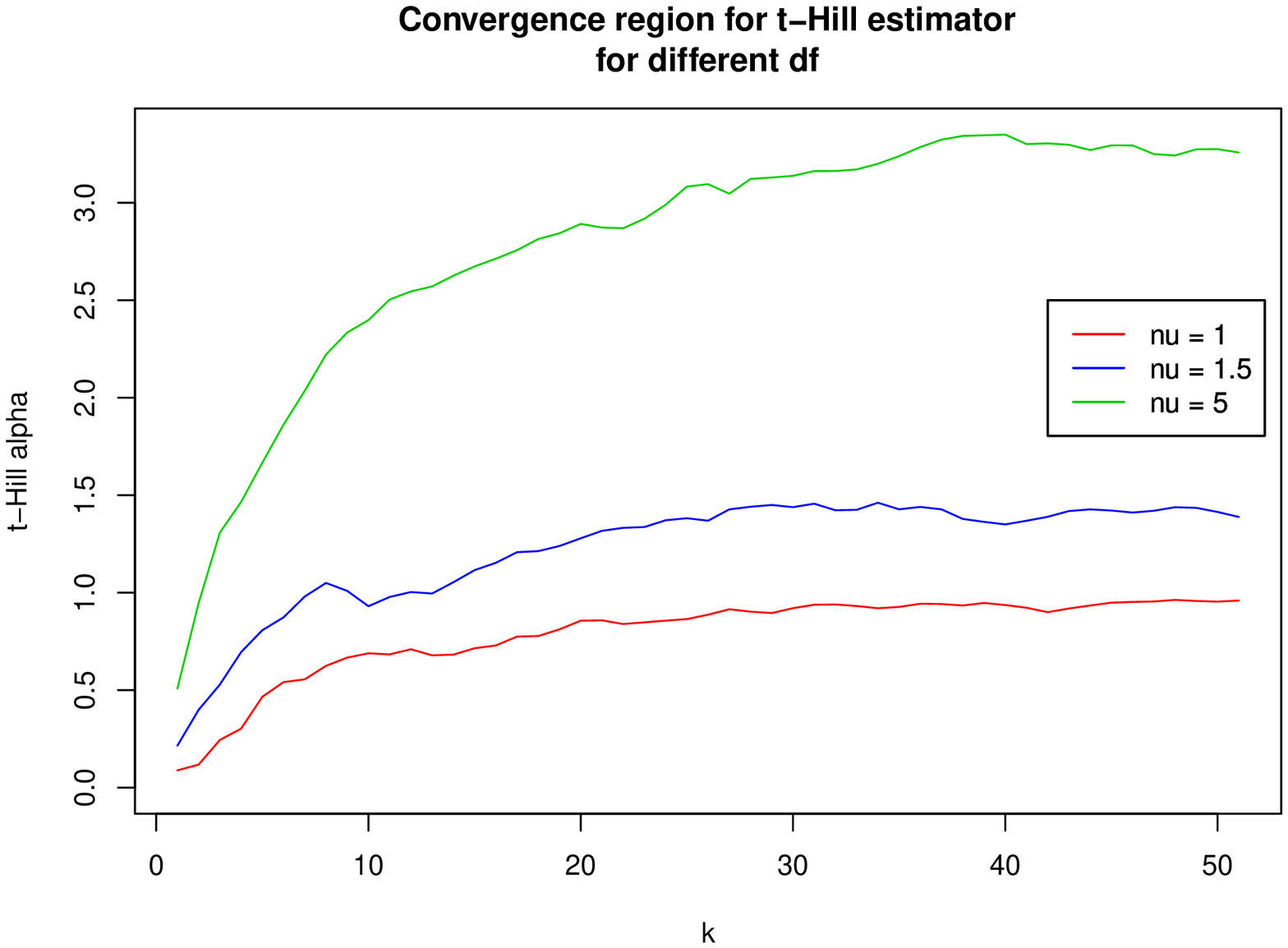}}
	\caption{Comparison and convergence region of 3 regimes of t-Hill plot}
	\label{fig:IGMMconv}
\end{figure}

The three colored areas representing Regimes are displayed in Figure \ref{fig:hillEstimatorColored}.
Therein also the reciprocal of the harmonic moment tail index estimator $H_{k, n}^*(1.001)$ (almost Hill-estimator) for unskewed LATAM data (by \texttt{IGMM} and \texttt{get.input} of \cite{Georg11}) is plotted as an estimate for $\hat\alpha$ (see \cite{Beran}). Simulations have shown that using unskewed or original LATAM data yields approximately the same results (not provided here), however, they differ slightly in the upper bound of k due to the occurrences of zeros leading to infinite values in the computations (division by zero). Therefore, we have provided the result for the skewed data, indeed, one way to correct for this obstacle can be to replace zeros by simulated values of uniform distribution between zero and the following order statistics of the returns, which is unequal to zero. Another way would be to use only values unequal to zero in order to avoid this problem.

A sample of length equal to the number of observations of the aforementioned data (n = 1421) has been simulated from Student distribution with $\nu = 5, 2, 1$ degrees of freedom. Following to that the harmonic moment tail index estimator has been computed on the basis of the ordered absolute values by setting $\beta = 2$ in Definition 1 on page 196 in  \cite{Beran}. This has been conducted for each degree of freedom in order to define the area of each Regime. These steps were repeated 100 times (Figure \ref{fig:IGMMconv} shows the result for 10 repetitions for the sake of comparison) for every setup and the reciprocal of the averages of these Hill estimators (in order to receive $\hat\alpha$) were plotted against the values of k, whereby $1 \leq k \leq n-1$. Recall that the almost Hill-estimator visualizes a single run of the algorithm, because it is based on the transformed absolute values of LATAM instead of simulated data.

It is well visible, that LATAM data tail is substantially overlapping with Regime III, thus it is not recommended to process these data with  \texttt{IGMM}.
To explain this fact from the point of view of finance,
we shall realize that LATAM data contains daily log-returns (in percent) of an equity fund investing in Latin America (LATAM) from January 1, 2002 until May 31, 2007. Emerging markets in Latin America
(see e.g. \cite{RamazanFaruk}) can have different properties on left and right tails. It was shown that e.g. Argentina and Brazil have higher estimates of the right tail index than of the left tail index. Therefore, high positive returns are more likely than similar losses in these growing economies. In 2004 it was observed that positive stock return distribution in e.g. Brazil may not have a finite second moment since the estimated extreme value index was around 0.5.
There is (even from 1988) an empirical evidence of non-existence of first moment (see page 9 of \cite{Akgiray}).
Another increase of heavy tailedness of the right tail has been introduced in the years 2004-2015, where high-frequency trading became more present in Latin America.
Analogously, in electric markets less credibility has been given to analysis using empirical means, often quickly replaced by median based techniques (see e.g. \cite{Kim11}).

\begin{figure}[h]
	\centering
		\includegraphics[width=0.8\textwidth]{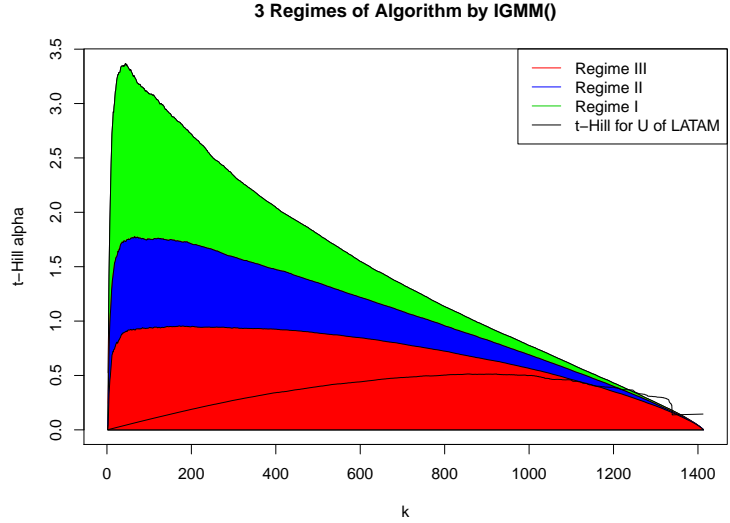}
	\caption{3 Regimes of t-Hill plot and \texttt{LATAM} data}
	\label{fig:hillEstimatorColored}
\end{figure}

\subsection{On Regime III of \texttt{IGMM}}

As mentioned above, Normal score is working in Regime I, but not in Regime III.
The classical score function as an indicator of the sensitivity of likelihood $L$, $S_\mu(X)= \frac{\partial}{\partial\theta} \log L(\theta;X),$
 has been built for distributions with support on real line, having all moments (see Fisher  (1925) \cite{Fisher}).
In case of Regime 3 (no finite moments), we shall not only transform a random variable, but also appropriately transform  its inference function.
For classical transformed t-score results see Fabi\'an (2001) \cite{Fabian2001} and Stehl\'\i k et al (2010) \cite{Stehlik08CS}.
In this letter we consider  only  {a} semi-parametric setup.
For a nonparametric analogy  see Dobrovidov, Koshkin and  Vasiliev (2012) \cite{Dobrovidov} where
scores $S_n=\frac{B^2}{A}\frac{\partial }{\partial x_n} \log f(x_n|x_{n-1})+\frac{x_n}{A}$ are defined for {a} conditionally exponential family
in the linear model $X_n = AS_n + B\eta_n,$ where $A,B$  are known constants, $\eta_n$ is Gaussian noise, $(X_n, S_n),n>1$ is {a} two-component Markov process,
$(X_n)$ is an observable process and $(S_n)$ is {an} unobservable useful process.

In our setup, let  ${\cal X}$ be {the} support of {the} distribution $F$ with density $f$,
continuously differentiable according to $x \in {\cal X}$ and let
$\eta: {\cal X} \to {\Re}$ be given by Johnson (1949) \cite{Johnson}
$\eta(x) =x,\ {\rm if}  {\cal X}={\Re}, \eta(x) =\log (x-a),\ {\rm if}  {\cal X}=(a,\infty)$ and $\eta(x)= \log \frac{x}{1-x},  \ {\rm if}  {\cal X}=(0,1).$
 Then the {\it
transformation-based score} or shortly the {\it t-score} (see Fabi\'an (2001) \cite{Fabian2001}) is defined
by $$T(x)=-\frac{1}{f(x)} \frac{d}{dx}\left(\frac{1}
 {\eta'(x)}f(x) \right), \label{T} $$
which expresses a relative change of a  "basic component of the density",  i.e., density divided by the Jacobian of  mapping $\eta.$

It is clear that for Normal distribution, which is an archetypical distribution, we have $\eta(x)=x, \ S(x, \theta) = \frac{d}{d \theta} \log f(x,\theta) $
and {$\hat{\theta} = {\rm MLE}$, with MLE standing for {\it maximum likelihood estimator,} which is the} solution of $ \sum_{i=1}^n S(X_i, \hat{\theta}) = 0. $

However, for {the} Pareto distribution we can consider two recently implemented approaches, namely:
\begin{itemize}

\item MLE, which {is related to  the ``standard score" } estimation with $\eta(x)=x$ and
$ S(X,\alpha) = \frac{1}{\alpha} - \log x $
 and

\item {$t$}-score estimation  with $\eta(x)=\log (x-1)$ (see Stehl\'\i k M. et al. (2010) \cite{Stehlik08CS}).
Notice that {the} MLE is not robust {wrt}  right outliers, i.e.\ {if} $ X_i \rightarrow \infty $, then $ \hat{\alpha} \downarrow 0.$
For {  $t$}-estimation we have \emph{t}-score $$ T(x) = \alpha \Big(1 - \frac{\alpha + 1}{\alpha x}\Big).$$
Thus standard estimation $ \sum T(X_i) = 0 $ {gives}  us
$ \hat{\alpha} = \frac{1}{\overline{x}-1} $ (where
$ \overline{x} = \frac{n}{\sum \frac{1}{x_i}} $ is harmonic mean) which is {an} estimator apparently robust against right-outliers.

\end{itemize}

Thus transformation of the data (e.g. by machinery of Lambert W variable), accompanied with a construction of proper score function transformation
is the  reasonable further research direction to regularize Algorithm 3 in Regime III.

\section{Skewness: On asset Returns and t-distribution}

Skewness and symmetry are fundamental objects of statistics and
it is interesting to study their transformations.
Symmetry itself is related to the nature of the problem and its permutation invariance, and cannot be obtained just by a simple transformation.
Thus symmetry is one of the fundamental notions of nonparametric statistics and is fundament for typical value of Hartigan (1969) \cite{Hartigan},
studied in perspective of reflection groups in    Francis, Stehl\'\i{k} and  Wynn (2014) \cite{Francis}.

Goerg G. M. (2011) \cite{Georg11} defines a transformation $Y_\gamma=U\cdot\exp(\gamma\cdot U)\sigma+\mu$
where $Y_\gamma$ is skewed output and $U$ is symmetrical input.
It is true, that having a symmetric zero-mean $U,$ $\gamma\neq 0$ regulates the skewness.
However, the inverse problem is much more delicate, as is demonstrated by the following simulations.
In Section 7.2, "Asset returns", Goerg G. M. (2011) \cite{Georg11} used Kolmogorov-Smirnov (KS) test, and stated "\emph{As a KS test cannot reject a student t-distribution..}".
KS test implementation in R\cite{R}, (as function \texttt{ks.test()}) was also used in the function \texttt{ks.test.t()} which was introduced in Goerg G. M. (2011) \cite{Georg11} and in his package \texttt{LambertW}. However, parameters $\hat \tau_{MLE}$ are estimated and thus, classical KS test cannot be used. There exist  some more refined distribution theory for the KS test with estimated parameters (see Durbin (1973) \cite{Durbin}), but this is not implemented in \texttt{ks.test()}, used in the function \texttt{ks.test.t()}. The undesirable parameter dependence of such implementation can influence one of the  goals of the paper: having a symmetric $t$-distribution input $U$ and $Y_\gamma=(U\exp(\gamma U))\sigma_x+\mu_x$ being a skewed output.

The following example shows, that estimation of parameters affects this aim in an undesirable way. First, we  simulated input variable $U$ as a skewed \emph{t}-distribution (see Fernandez and Steel (1998) \cite{FernandezSteel1998}) with skew parameter $\gamma^*$.  Data was simulated with function \texttt{rskt(n, df, $gamma^*$)} of package \texttt{skewt}. The values for parameter $\gamma^*$ and resulting skewness with four degrees of freedom can be found in left part of Table \ref{tab:skewp}. Then we transformed data to $Y=(U\exp(-bU))c+a,$ where $a,b$ and $c$ have been chosen from grids $a = seq(0,1,by = 0.01); b = seq(0,1,by = 0.01); c = seq(0.1,1.5,by = 0.01).$ Finally we  estimated $U$ and parameters by \texttt{IGMM} and conducted \texttt{ks.test.t()} from this package \texttt{LambertW}. This shows the effect of usage of \texttt{ks.test.t()} jointly with parameter estimation, which led to acceptance of skewed distributions as symmetric student distribution.

For $\gamma^*$ equal to 0.9 or 0.75 we received p-values of 0.502 and 0.269 , and thus skewed distribution (skewness = -0.93 and -1.37) is accepted as symmetric student. First line of Table \ref{tab:skewp} presents simulation of \emph{t}-distribution (skewness = -0.3415) and resulting p-value is correctly higher than 0.05. The same comparison was done for skewed normal distribution, which was simulated with \texttt{rsn(n, xi , omega , alpha )} from package \texttt{sn}. $\alpha$ is in this setting skewing parameter and its values can be seen in the first column of the right side of Table \ref{tab:skewp}. Location (\texttt{xi}) and scaling parameters (\texttt{omega}), which are equivalent to mean and standard deviation, were chosen to be $\mu = 4$ and $\sigma = 2$. Skewness was compared for different $\alpha$ and p-values resulting from \texttt{ks.test.t()}  are presented as before. For all listed cases we received p-values $p > 0.05$  and therefore skewed normal distributions were falsely assumed as symmetric \emph{t}-distributions.  Obviously p-values are decreasing for higher $\alpha$, but for all $\alpha \leq 8$  symmetric \emph{t}-distribution cannot be rejected for simulated skewed normal distribution.

\begin{table}[h]
\begin{center}
\caption{Skewness and p-values resulting from \texttt{ks.test.t()} of simulated t-, skewed t-, normal and skewed normal distribution. $\gamma^*$ (t) and $\alpha$ (normal) are skewing parameters. Parameters of un-skewed distributions are in first row.}
\label{tab:skewp}
\begin{tabular}{p{1.32cm}p{1.32cm}p{1.32cm}|p{1.79cm}p{1.79cm}p{1.79cm}}\hline
\multicolumn{3}{c|}{\textbf{Skewed t- and t-distribution}} & \multicolumn{3}{c}{\textbf{Skewed normal and normal distribution}} \\
 skewness & p-value & $\gamma^*$ & skewness & p-value & $\alpha$ \\
  \hline
 -0.3415 & 0.2872 &   & 0.0130 & 0.7731 & \\
 -2.8894 & 0.0001 & 0.20 &  -0.0892 & 0.0569 & 0.10 \\
 -1.6850 & 0.0000 & 0.40 &  0.0277 & 0.5323 & 0.50 \\
 -1.3785 & 0.2693 & 0.75 & 0.0810 & 0.6801 & 1.00 \\
  -0.9304 & 0.5017 & 0.90 &0.8108 & 0.5035 & 2.50 \\
& & &  0.8054 & 0.0924 & 5.00 \\
& & &  0.9391 & 0.0527 & 8.00 \\
   \hline
\end{tabular}
\end{center}
\end{table}

In order to check graphically for impact of $\gamma^*$ and $\alpha$ on skewed distributions, kernel density estimations were plotted in R \cite{R}. This density estimation comparison in Figure \ref{fig1}(a) shows stronger skewed distributions for decreasing values of $\gamma^*$. Distributions were simulated with negative skewness in this example. Black density corresponds to student distribution and the others are computed with previously defined $\gamma^*$ values and show skewed \emph{t}-distributions. A graphical comparison between skewed normal distributions and normal distribution is done in Figure \ref{fig1}(b). Increasing skewing parameter $\alpha$ leads to stronger skewness of the data and a shift to the right. In contrast to the previous examples skewness is except for $\alpha = 0.1$ positive and increasing with $\alpha$.

\begin{figure}[h]
%\vspace{6pc}
		\subfigure[Skewed \emph{t}- and \emph{t}-distribution, four degrees of freedom]{\includegraphics[width = 0.48\textwidth]{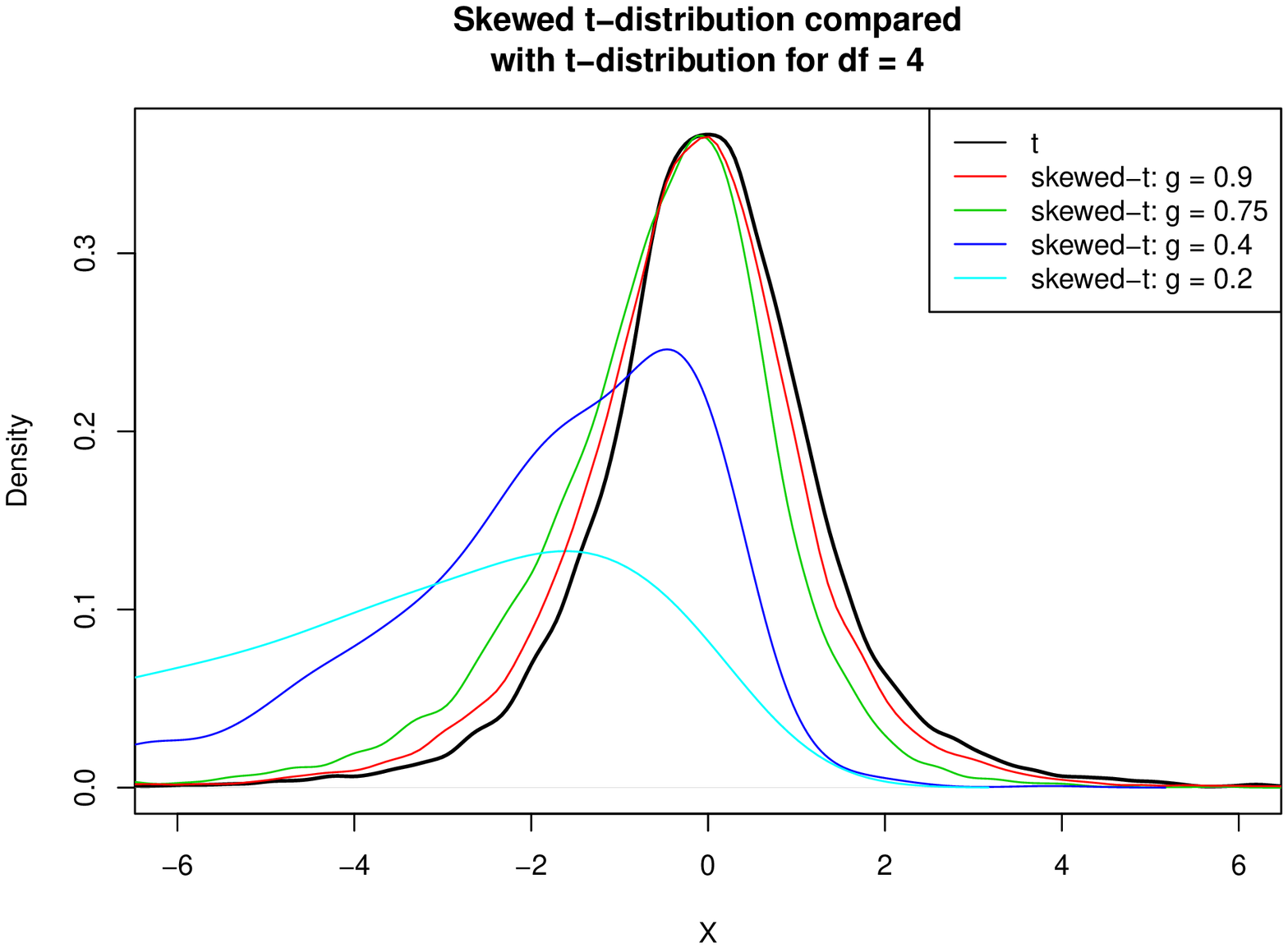}}
		\subfigure[Skewed normal and normal-distribution, $\mu = 4$ and $\sigma = 2$]{\includegraphics[width = 0.48\textwidth]{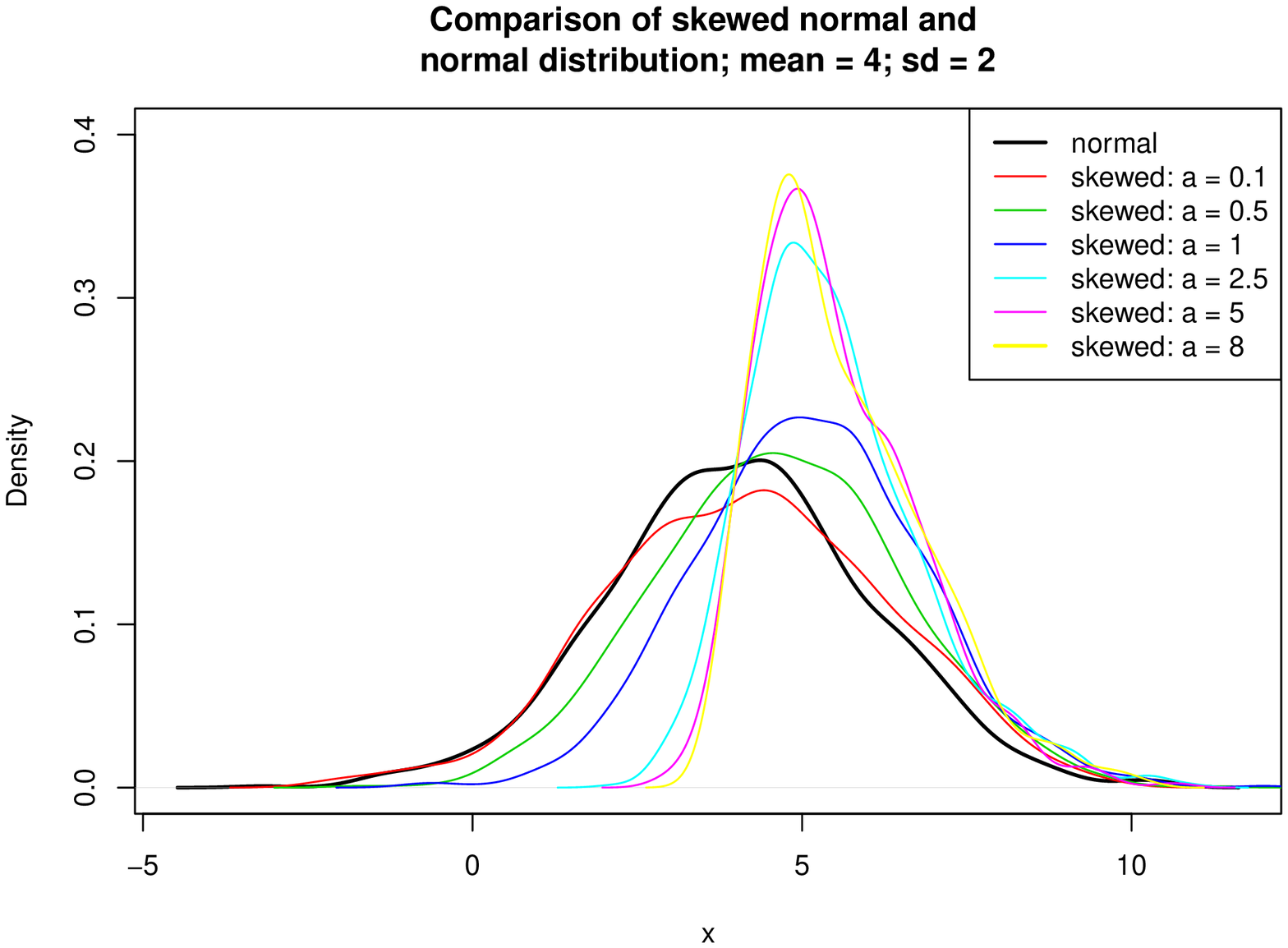}}
	\caption{Comparison of skewed and unskewed t(df = 4)- and normal($\mu=4, \sigma=2$)-distribution}
	\label{fig1}
\end{figure}

\subsubsection{Auto-Correlation Rising from \texttt{IGMM} and LATAM data}
We also checked auto-correlations resulting from estimation by Algorithm 3 for different distributions in a following simulation setup. We simulated  standard Normal distribution, Weibull, Exponential and student-\emph{t} distributions. In the next step \texttt{IGMM} was used to estimate parameters and as a consequence back-transformation with \texttt{get.input} was applied with estimated $\hat\mu, \hat\sigma$ and $\hat\gamma$.
We observed significant auto-correlation for all 4 distributions.
Also  auto-correlation function of back-transformed series of LATAM has been observed to  be significant (e.g. at lags 2, 7, 8, 13 and 30).

\bibliographystyle{imsart-nameyear}

\begin{thebibliography}{9}

\bibitem[Akgiray et al. (2014)]{Akgiray} \textsc{Akgiray, V., Booth, G.G., Seifert, B.} (1988)    Distribution properties of Latin American black market exchange rates. \textit{Journal of International Money and Finance}, \textbf{7}(1),  37-48

\bibitem[Balakrishnan  and Stehl\'\i k (2008)]{Bala} \textsc{Balakrishnan N. and Stehl\'\i k M.} (2008), Exact likelihood ratio test of the scale for censored Weibull sample. \textit{Ifas Res.Report}  \textbf{35}, online at http://www.jku.at/ifas/

\bibitem[Beran et al. (2014)]{Beran}  \textsc{Beran J.,  Schell D.,   Stehl\'\i k M.} (2014)  The {harmonic moment tail index estimator}: asymptotic distribution and robustness, \textit{Annals of the Institute of Statistical Mathematics} \textbf{66}(1); 193-220.



\bibitem[Dobrovidov et al. (2012)]{Dobrovidov} \textsc{Dobrovidov A.V., Koshkin G.M.\ and  Vasiliev V.A} (2012) \textit{Non-Parametric State Space Models}. Kendrick Press. USA

\bibitem[Durbin, J. (1973)]{Durbin} \textsc{Durbin, J.} (1973) \textit{Distribution theory for tests based on the sample distribution function}. SIAM.

\bibitem[Fabi\'{a}n (2001)]{Fabian2001} \textsc{Fabi\'{a}n Z} (2001).  Induced cores and their use in robust
parametric estimation.  \textit{Communications in Statistics---Theory Methods} {\bf 30} 537--556

\bibitem[Fabi\'{a}n and  Stehl\'\i{k} (2009)]{FSt09} \textsc{Fabian, Z. and Stehl\'\i{k}, M.} (2009). On robust and distribution sensitive Hill like method. \textit{IFAS Research Paper Series} $\bf{43(4)}$, online at http://www.jku.at/ifas/

\bibitem[Fabi\'{a}n and  Stehl\'\i{k} (2009)]{BGP2012} \textsc{Brilhante, M.F., Gomes, M.I., Pestana, D.} (2013), A simple generalization of the Hill estimator, Computational Statistics \& Data Analysis {\bf 57}, 518--535

\bibitem[Fernandez, C. and Steel, M. F. J. (1998)]{FernandezSteel1998} \textsc{Fernandez, C. and Steel, M. F. J.} (1998). On Bayesian modeling of fat tails and skewness, \textit{J. Am. Statist. Assoc.} \textbf{93}  359-371.

\bibitem[Fisher  (1925)]{Fisher}  \textsc{Fisher R. A.} (1925). Theory of statistical estimation, \textit{Proceedings of the Cambridge Philosophical Society} \textbf{22} 700-725,  doi:10.1017/S0305004100009580

\bibitem[Francis et al.  (2014)]{Francis}   \textsc{Francis, R.A., Stehl\'\i{k}, M.  and  Wynn, H.P.} (2014) Exact confidence nets based on finite reflection groups, arXiv:1407.8375 [math.ST]

\bibitem[Gen\c cay and Sel\c cuk (2004)]{RamazanFaruk} \textsc{Gen\c cay, R. and Sel\c cuk, F.} (2004). Extreme value theory and Value-at-Risk: Relative performance in emerging markets. \textit{International Journal of Forecasting},  \textbf{20}(2), 287-303.


\bibitem[Goerg (2011)]{Georg11}  \textsc{Goerg G. M.} (2011).  Lambert W Random Variables - A New Family Of Generalized Skewed Distributions With Applications To Risk Estimation. \textit{The Annals of Applied Statistics}. \textbf{5}(3) 2197-2230.

\bibitem[Gomes et al. (2014)]{VARMOP}  \textsc{Gomes, M.I.,   Brilhante, F. and Pestana, D.}(2014). A mean-of-order-p class of value-at-risk estimators.  Theory and Practice of Risk Assessment, Springer Proceedings in Mathematics and Statistics, In Kitsos, C., Oliveira, T., Rigas, A. and Gulati, S. (eds.), p. 1-16,

\bibitem[Hartigan (1969)]{Hartigan} \textsc{Hartigan, J. A.}(1969)
Using subsample values as typical values. \emph{Journal of the
American Statistical Association} \textbf{64}, 328, 1303-1317.

\bibitem[Hill (1975)]{Hill1975} \textsc{Hill, B.}(1975)  A simple general approach to inference about the tail of a distribution, \textit{Annals of Statistics} \textbf{3}:5, 1163--1174



%\bibitem[Cizek et al. (2005)]{Heardle} \textsc{Cizek et al. }(2005) \textit{Statistical Tools for Finance and Insurance}, Eds. Cizek, P., H\"ardle, W. K., Weron, R. (Eds.) , Springer-Verlag.

\bibitem[Johnson (1949)]{Johnson} \textsc{Johnson N. L.} (1949). Systems of frequency curves generated by methods of
translations. \textit{Biometrika} \textbf{36} 149-176.

\bibitem[Jordanova et al. (2013)]{Chemo1} \textsc{Jordanova P., Du\v sek J. and Stehl\'\i k M.} (2013), Modeling methane emission by the infinite moving average process, \emph{Chemometrics and Intelligent Laboratory Systems}, \textbf{ 122}, 40-49

\bibitem[Kim et al. (2011)]{Kim11} \textsc{Kim, J. H., Powell, W. B., and Collado, R. A.} (2011). Quantile optimization for heavy-tailed distribution using asymmetric signum functions. \textit{Princeton University}.

\bibitem[Lambert (1758)]{Lambert} \textsc{Lambert JH}(1758). \textit{Observationes variae in mathesin puram}. Acta Helveticae physico-mathematico-anatomico-botanico-medica, Band III, 128--168.



\bibitem[P\'{a}zman (1993)]{Pazman1} \textsc{P\'{a}zman  A,}(1993). {\textit Nonlinear statistical Models.} Kluwer Acad. Publ. Dordrecht. chapters 9.1 and 9.2

\bibitem[R (2008)]{R} \textsc{R Core Development Team} (2008): A Language and Environment for Statistical Computing. R Foundation for Statistical Computing. Vienna, Austria. ISBN 3-900051-07-0

\bibitem[{Stehl\'\i{k}(2003)}]{Stehlik03} \textsc{Stehl\'\i{k}, M.} (2003).  Distributions of exact tests in the exponential family. \textit{Metrika} \textbf{57} 145--164.

\bibitem[Stehl\'\i k (2006)]{Stehlik06} \textsc{Stehl\'\i k M.} (2006). Exact likelihood ratio scale and homogeneity testing of some loss processes. \textit{Statistics and Probability Letters} \textbf{76} 19-26.

\bibitem[{Stehl\'\i{k}(2008)}]{Stehlik08} \textsc{Stehl\'\i{k}, M.} (2008). Homogeneity and scale testing of generalized gamma distribution.  \textit{Reliability Engineering \& System Safety} \textbf{93} 1809--1813.

\bibitem[{Stehl\'\i{k} et al. (2010)}]{Stehlik08CS} \textsc{Stehl\'\i k M., Potock\'y R.,   Waldl, H. and   Fabian, Z.} (2010). On the
favourable estimation of fitting heavy tailed data. \textit{Computational Statistics} \textbf{25} 485-503

\bibitem[{Stehl\'\i{k} et al. (2012) }]{StehlikCS} \textsc{Stehl\'\i k M., Fabi\'{a}n Z. and  St\v relec L.} (2012). \textit{Small sample robust testing for Normality against Pareto tails}.  \textit{Communications in Statistics - Simulation and Computation} \textbf{41}(7) 1167-1194

\bibitem[{Stehl\'\i{k} et al. (2014) }]{StehlikCIS}  \textsc{Stehl\'\i k, M. Thulin, M.,  St\v relec, L.} (2014). On robust testing for normality in chemometrics. \textit{Chemometrics and Intelligent Laboratory Systems} \textbf{130}  98-108

\bibitem[Sung (2013)]{Sung}  \textsc{Sung  S.H.}  (2013). On the strong law of large numbers for pairwise i.i.d. random
variables with general moment conditions,  \emph{Statistics and Probability Letters} \textbf{83} 1963-1968

\bibitem[Wilks  (1938)]{Wilks}  \textsc{ Wilks S. S.} (1938). The Large-Sample Distribution of the Likelihood Ratio for Testing Composite Hypotheses. \textit{Ann. Math. Statist.} \textbf{9}(1) 60-62.

\end{thebibliography}

\end{document}